\documentclass[11pt]{amsart}

\usepackage[T2C]{fontenc}
\usepackage[utf8]{inputenc}
\usepackage[english]{babel}
\usepackage{amsfonts,amssymb,amscd}
\usepackage{mathrsfs}
\usepackage{graphics}

\usepackage{hypbmsec}
\usepackage{float}

\usepackage{longtable}
\usepackage{booktabs}
\usepackage[hidelinks]{hyperref}
\usepackage[numbers]{natbib}

\newtheorem{Theorem}{Theorem}

\newtheorem{Remark}{Remark}

\newtheorem{Example}{Example}

\newcommand{\eps}{\varepsilon}
\newcommand{\leg}[1]{\left(\dfrac{#1}{p}\right)}
\DeclareMathOperator{\End}{End}

\title{Quadratic residue patterns of length 4 and 5}

\author{Abdulkadyr Buchaev}
\address{
Higher School of Modern Mathematics MIPT,
Moscow, Russia
}
\email{buchaev.aia@phystech.edu}

\author{Michael Tsfasman}
\address{
Higher School of Modern Mathematics MIPT,
Moscow, Russia
}
\email{mtsfasman@yandex.ru}


\thanks{The research is supported by MSHE RF GZ project}

\begin{document}

\begin{abstract}
We generalize the results of the first part of the paper by Kiritchenko, Tsfasman, Vl\u{a}du\c{t} and Zakharevich on quadratic residue patterns to the case of arbitrary words over the alphabet $\{R, N\}$ of length $l = 4$ and $l = 5$. First, we obtain explicit formulas for the number of occurrences of a given pattern $S$ in the sequence of quadratic residues and nonresidues modulo a prime $p$. These formulas are expressed in terms of Frobenius traces on a finite collection of low-genus hyperelliptic curves. Second, using a generalized Sato--Tate approach, we describe the limiting distributions of the normalized error term for $l = 4$.
\end{abstract}

\maketitle

\section{Introduction}

Let $p$ be an odd prime. Since the time of Fermat, it has been known that among the nonzero residue classes $1, \ldots, p-1$ modulo $p$, the number of quadratic residues equals the number of quadratic nonresidues. Squares are referred to as (quadratic) residues, and denoted by $R$, whereas nonsquares as (quadratic) nonresidues, and denoted by $N$. At the end of the nineteenth century N. Aladov \cite{l2aladov} asked for, and determined, the numbers of of occurences of the four pairs $(RR), (RN), (NR)$ and $(NN)$. The same problem concerning words in $R$ and $N$ of arbitrary fixed length $l$ was solved for $l=3$ by E. Jacobsthal in his doctoral thesis. For a review of subsequent work, we refer the reader to the paper by Kiritchenko, Tsfasman, Vl\u{a}du\c{t} and Zakharevich \cite{main_paper}.  In \cite{main_paper}, the authors treated only the pattern $(RRRR)$ of four consecutive quadratic residues. The main novelty of that paper was its substantial use of algebraic geometry.

For an odd prime $p$, we form a word in $R$ and $N$ by replacing each integer in $(1, \ldots, p - 1)$ with $R$ if it is a quadratic residue and with $N$ otherwise. 
\begin{equation}
W_p = (w_p(1), w_p(2), \ldots, w_p(p-1)) \in \{R, N\}^{p-1},
\end{equation}
\begin{equation}
w_p(n) = \begin{cases} R, \quad \left(\frac{n}{p}\right) = 1, \\
					N, \quad \left(\frac{n}{p}\right) = -1,
\end{cases}
\end{equation}
where $\left(\frac{n}{p}\right)$ denotes the Legendre symbol of an integer $n$ modulo $p$.

For a fixed pattern (i.e., a word of length $l$) $S = s_1\ldots s_l \in \{R, N\}^l$ define the number
\begin{equation}
n_p(S) = \#\{1 \le x \le p - l : (w_p(x), \ldots, w_p(x+l-1)) = S\}
\end{equation}
of subwords of the word $W_p$ that are equal to $S$.

\begin{Example}
Let $p = 13$, $l = 3$. Then $W_p = (RNRRNNNNRRNR)$,
\begin{equation}
n_p(RRR) = 0, n_p(RRN) = 2, n_p(RNR) = 2, n_p(NRR) = 2,
\end{equation}
\begin{equation}
n_p(NNR) = 1, n_p(NRN) = 0, n_p(RNN) = 1, n_p(NNN) = 2.
\end{equation}
\end{Example}

How can $n_p(S)$ be calculated? For $l = 1$ one it is known that $n_p(R) = n_p(N) = \frac{p-1}{2}$. The case $l = 2$ was first considered by Aladov in his paper \cite{l2aladov}, where he gets
\begin{equation}
n_p(RR) = \dfrac{p - 4 - \leg{-1}}{4}, n_p(RN) = \dfrac{p - \leg{-1}}{4},
\end{equation}
\begin{equation}
n_p(NN) = n_p(NR) = \dfrac{p - 2 + \leg{-1}}{4}.
\end{equation}

In \cite{l3jacobsthal}, E. Jacobsthal solved the problem for $l=3$. Namely, for $p = 4k+3$ he gets
\begin{equation}
n_p(RRR) = n_p(NNN) = n_p(NRR) = n_p(NNR) = \dfrac{p-5-2\leg{2}}{8},
\end{equation}
\begin{equation}
n_p(NRN) = n_p(RNR) = n_p(RNN) = n_p(RRN) = \dfrac{p-1+2\leg{2}}{8}.
\end{equation}
The formulas for $p = 4k+1$ are more complicated. Denote
\begin{equation}
a = \sum_{i = 1}^{p-3} \leg{i(i+1)(i+2)}.
\end{equation}
Then he gets
\begin{equation}
n_p(RNR) = n_p(RRN) = n_p(NRR) = n_p(NNN) = \dfrac{p-3-a}{8},
\end{equation}
\begin{equation}
n_p(RNN) = n_p(NNR) = \dfrac{p+1+a}{8},
\end{equation}
\begin{equation}
n_p(NRN) = \dfrac{p-3+4\leg{2}+a}{8},
\end{equation}
\begin{equation}
n_p(RRR) = \dfrac{p-11-4\leg{2}+a}{8}.
\end{equation}

As we shall see below, the problem is connected with counting $\mathbb{F}_p$-points on certain curves. For $l=2$ the corresponding curve is rational, and thus Aladov's answer is so simple. For $l=3$ the curve is elliptic, which explains why Jacobsthal's answer is much more complicated. 

The research was continued, and in the paper \cite{main_paper} the authors consider the general situation and establish a geometric framework that explains formulas that occur for any $l$. In particular, they  study in detail the case of the \textit{homogeneous} patterns $(RRRR)$ of length $l = 4$. This is much more difficult, since the genus of the corresponding curve is $5$. However, the answer for $l = 4$ can be obtained due to some specific properties of this curve.

In particular, they get

\begin{equation}\label{p30}n_p(RRRR)=\frac{p}{16}-\frac{a_4(p)}{16} +c_p(4) \;\mbox{for}\; p=4k+3, \end{equation}
\begin{equation}\label{p10}n_p(RRRR)=\frac{p}{16}-\frac{a_0(p)}8-\frac{a_1(p)}8-\frac{a_4(p)}{16}+c_p(4) \;\mbox{for}\; p=4k+1, \end{equation}
where $a_i(p)$ are the Frobenius traces of certain elliptic curves $E_i$ over $\mathbb{F}_p$; in particular, $a_0(p)$ is the quantity $a$ of Jacobsthal for $(RRR)$, and $a_4$ is its direct generalization
\begin{equation}
a = \sum_{i = 1}^{p-3} \leg{i(i+1)(i+2)(i+3)}.
\end{equation}

Here the function $c_p(4)$ of $p$ is bounded by an absolute constant $c(4)$. 
We compute it explicitly in Theorem 1; it depends only on $p \;\mbox{mod}\;24$.

Just as for $l=3$, the case of $p=4k+1$ is more complicated than that of $p=4k+3$.

At the end of the first part of \cite{main_paper} and in \cite{vladuc} the authors also investigate statistical properties of the deviation from the linear term $p/{2^l}$ of $n_p(S)$ under the \textit{generalized Sato-Tate conjecture} (GST). Our aim is to generalize both the geometric and statistical results of these works to the case of \textit{arbitrary} patterns of lengths $4$ and $5$.

The first part of our work gives explicit formulas for $n_p(S)$ for arbitrary patterns $S$ of lengths $4$ and $5$ in terms of Frobenius traces and explicit correction terms coming from zeros of the Legendre symbol. Our second main result describes the limit distributions of the corresponding normalized error terms for $l = 4$ via the generalized Sato--Tate framework.

This paper is organized as follows. In Section 2 we present explicit formulas for counting $n_p(S)$ in terms of Frobenius traces and correction terms. In Sections 3 and 4 we list elliptic and hyperelliptic curves whose traces occur in formulas for $l = 4$ and $l = 5$, respectively, and establish relations between these curves. In Section 5 we describe the limit distribution of the normalized error terms for $l = 4$ using GST.

\section{Explicit formulas}

In this section we will write down explicit formulas for $n_p(S)$.
These formulas can be found in \cite{cohen_triples_quadruples} and \cite{conrad_quadratic_residue_patterns}, but we will write them down for the convenience of the reader. 
Choose a fixed  quadratic residue and nonresidue $\eps_+$ and $\eps_-$, respectively, in $\mathbb{F}_p$. We then identify $S$ with an element of $\{\eps_+, \eps_-\}^l$ in the following way:
\begin{equation}
S = \eps_1\ldots \eps_l, \qquad\eps_i = \begin{cases}
								\eps_+, \quad s_i = R, \\
								\eps_-, \quad s_i = N.
\end{cases}
\end{equation}

Define a function $I_S(x): \mathbb{F}_p \to 2^{-l}\mathbb{Z}$ as follows:
\begin{equation}
I_S(x) = \prod_{i = 1}^l \frac12 \left(1 + \left(\dfrac{\eps_i \cdot (x + i - 1)}{p}\right)\right).
\end{equation}
It is easy to see that for $1 \le x \le p - l$ the function $I_S(x)$ is the indicator of the event that the subword $(w_p(x), \ldots, w_p(x + l -1))$ equals $S$. Then we have
\begin{equation}
n_p(S) = \sum_{x = 1}^{p - l} I_S(x).
\end{equation}
Now, let us rewrite this sum as follows:
\begin{equation}
n_p(S) = \sum_{x = 1}^{p - l} I_S(x) = \sum_{x \in \mathbb{F}_p} I_S(x) - \sum_{x = p - l + 1}^{p} I_S(x) =
\end{equation}
\begin{equation}
= 2^{-l} \left(\sum_{x \in \mathbb{F}_p} \prod_{i =1}^l \left(1 + \leg{\eps_i \cdot (x + i -1)}\right)\right) - c_{p, S}(l) =
\end{equation}
\begin{equation}
 = 2^{-l} \left(\sum_{x \in \mathbb{F}_p} \left(1 + \sum_{t = 1}^l \sum_{\substack{T \subset [1, l] \\ |T| = t}} \prod_{i \in T} \leg{\eps_i \cdot (x + i - 1)}\right)\right) - c_{p, S}(l) =
\end{equation}
\begin{equation}
= 2^{-l} \left(p + \sum_{t = 1}^l \sum_{\substack{T \subset [1, l] \\ |T| = t}} \sum_{x \in \mathbb{F}_p}\prod_{i \in T} \leg{\eps_i \cdot (x + i -1)}\right) - c_{p, S}(l) =
\end{equation}
\begin{equation}\label{concise_formula}
= 2^{-l}p + \sum_{t = 1}^l \sum_{\substack{T \subset [1, l] \\ |T| = t}} N_{T, S} - c_{p, S}(l),
\end{equation}
where 
\begin{equation}
N_{T, S} = \sum_{x \in \mathbb{F}_p} \prod_{i \in T} \leg{\eps_i \cdot (x + i - 1)} = \sum_{x \in \mathbb{F}_p} \leg{f_{T, S}(x)} \in \mathbb{Z}
\end{equation}
for
\begin{equation}
f_{T, S}(x) = \prod_{i \in T} \varepsilon_i (x + i) \in \mathbb{Z}[x],
\end{equation}
and
\begin{equation}
c_{p, S}(l) = \sum_{x = p - l +1}^p I_S(x) = \sum_{x = p-l+1}^p \prod_{i = 1}^l \frac12 \left(1 + \leg{\eps_i\cdot (x + i -1)}\right).
\end{equation}

It follows that $N_{T, S}$ is closely related to the number of points over $\mathbb{F}_p$ of the affine curve $C^a_{S, T}$ given by the equation $y^2 = f_{S, T}(x)$. Namely, if $|C(\mathbb{F}_p)|$ is the number of points of a curve $C$ over $\mathbb{F}_p$, then $N_{S, T} = |{C_{S, T}^a}(\mathbb{F}_p)| - p$.

Let $C_{S, T}$ be the smooth projective model of $C_{S, T}^a$. The genus of ${C_{S, T}}$ is equal to $\lfloor (t-1)/2\rfloor$, whence $C_{S, T}$ is elliptic for $t = 3, 4$ and hyperelliptic for $t \ge 5$.

We have
\begin{equation}
|{C_{S, T}}(\mathbb{F}_p)| = |C^a_{S, T}(\mathbb{F}_p)| + r_{S, T},
\end{equation}
 where
\begin{equation}
r_{S, T} = \begin{cases}
1, & |T| \text{ is odd,} \\
1 + \leg{\prod_{i \in T} \eps_i}, & |T| \text{ is even.}
			\end{cases}
\end{equation}

We also have $|{C_{S, T}}(\mathbb{F}_p)| = p + 1 - a_{S, T}$, where $a_{S, T}$ is the \textit{Frobenius trace} of the curve ${C_{S, T}}$, i.e., it is the trace of the Frobenius operator acting on the \'etale $\ell$-adic cohomology group $H^1({C_{S, T}}, \mathbb{Q}_\ell)$ or its inverse acting on the Tate module $T_\ell({J_{S, T}})$ of its Jacobian ${J_{S, T}}$. 

Therefore, we may rewrite equation (\ref{concise_formula}) as
\begin{equation}
n_p(S) = 2^{-l}p + 2^{-l} \sum_{\varnothing \ne T \subset [1, l]} (1 - a_{S, T} - r_{S, T}) - c_{p, S}(l) =
\end{equation}
\begin{equation}
 = 2^{-l}p - 2^{-l}\sum_{\varnothing \ne T \subset [1, l]} a_{S, T} + 2^{-l} \sum_{\varnothing \ne T \subset [1, l]}(1 - r_{S, T}) - c_{p, S}(l).
\end{equation}

Note that $a_{S, T}$ is zero for $t<3$, since $C_{S, T}$ is a rational curve.

Thus, the problem can be viewed in terms of sums of Legendre symbols as in earlier works or in terms of sums of traces of Frobenius. The latter interpretation arises from the following construction.

Fix a prime $p$, a quadratic residue $\eps_+$ and a quadratic nonresidue $\eps_-$ modulo $p$. For a pattern $S = \eps_1 \ldots \eps_l$ consider the affine curve $C_S^a$ over $\overline{\mathbb{F}_p}$ defined as the intersection of $l - 1$ quadrics:
\begin{equation}
\eps_2 x_2^2 - \eps_1 x_1^2 = 1, \quad \ldots, \quad \eps_l x_l^2 - \eps_{l-1} x_{l-1}^2 = 1. 
\end{equation}
Every point $(x_1, \ldots, x_l)$ of the curve $C_S^a$ yields the $l$-tuple $(\eps_1 x_1^2, \ldots, \eps_l x_l^2)$ that forms the pattern $S$. Since $x_i$ and $-x_i$ produce the same term, we obtain
\begin{equation}\label{number_of_points_and_np}
n_p(S) = 2^{-l}|C_S^{a,\circ}(\mathbb{F}_p)|,
\end{equation}
where $C_S^{a,\circ} \subset C_S^a$ consists of points whose coordinates are nonzero. Hence, the curve $C_S^a$ controls $n_p(S)$.

Let $C_S$ be the smooth projective model of $C_S^a$. In \cite[Theorem 2.3]{main_paper} the authors proved that the Jacobian of $C_{R \ldots R}$ is isogenous to the product of elliptic and hyperelliptic Jacobians of the curves ${C_{R \ldots R, T}}$, $|T|\ge 3$. Note that the proof can be adapted to any pattern $S$, not necessarily $S = (R \ldots R)$. This shows that the trace of Frobenius of $C_S$ is equal to the sum of traces of Frobenius of ${C_{S, T}}$, which explains the formulas obtained above. From the adjunction formula we obtain that the genus of $C_S$ is equal to $2^{l-2}(l-3) + 1$.

In the following sections we will simplify the error term $\sum_{\varnothing \ne T \subset [1, l]} a_{S, T}$ and present explicit values of the correction term $2^{-l}\sum_{\varnothing \ne T \subset [1, l]} (1 - r_{S, T}) - c_{p, S}(l)$ depending on $p\;\mbox{mod}\;24$ for all patterns $S$ of lengths $4$ and $5$, respectively.

\section{Case $l = 4$}

In this section we establish relations between the curves $C_{S, T}$ that allow us to simplify error terms $\sum_{\varnothing \ne T \subset [1, l]} a_{S, T}$. Simplification is due to the fact that curves that are isomorphic over $\mathbb{F}_p$ have the same Frobenius trace, whereas the Frobenius traces of a curve and its quadratic twist differ by a sign.

Note that for $|T| \le 2$ all curves $C_{S, T}$ are rational, hence $a_{S, T} = 0$. For $|T| \ge 3$ consider the following five curves $E_i$ that are smooth projective closures of the affine curves given by the following equations:
\begin{equation}
E_0: y^2 = x(x+1)(x+2); \quad E_1: y^2 = x(x+1)(x+3);
\end{equation}
\begin{equation}
E_2: y^2 = x(x+2)(x+3); \quad E_3: y^2 = (x+1)(x+2)(x+3); 
\end{equation}
\begin{equation}
E_4 : y^2 = x(x+1)(x+2)(x+3).
\end{equation}

Recall that $C_{S, T}^a$ is given by the equation $y^2 = f_{S, T}(x) = \prod_{i \in T} \eps_i (x+i)$. The curve $C_{S, T}^a$ is therefore isomorphic to the curve given by the equation $\prod_{i \in T} \eps_i y^2 = \prod_{i \in T} (x+i)$. Hence, $C_{S, T}$ is isomorphic over $\mathbb{F}_p$ to some $E_j$ if $\prod_{i \in T} \eps_i$ is a quadratic residue modulo $p$ or isomorphic to the quadratic twist of some $E_j$ otherwise. There is a one-to-one correspondence between the set  $\{E_i\}$ and the set $\{C_{S,T} \}$ which allows us to rewrite the error term $\sum_{\varnothing \ne T \subset [1, l]} a_{S, T}$ in terms of the traces $a_j$ of the curves $E_j$.

First, one have $a_0 = a_3$ since $E_0$ and $E_3$ are isomorphic over $\mathbb{Q}$; curves $E_1$ and $E_2$ are isomorphic over $\mathbb{F}_p$ if $p = 4k+1$ and $E_1$ is the quadratic twist of $E_2$ if $p = 4k+3$, whence $a_1 = (-1)^\frac{p-1}{2} a_2$. Also, $E_0$ is isomorphic to its own quadratic twist if $p = 4k+3$ which leads to $a_0 = 0$ in this case.

Combining these observations with the preceding discussion, we obtain the following theorem.
\begin{Theorem} For a prime $p \ge 5$ and a pattern $S$ of consecutive quadratic residues and quadratic nonresidues of length $l = 4$, the number $n_p(S)$ of occurrences of $S$ in $W_p$ is controlled (in the sense of equation (\ref{number_of_points_and_np})) by the curve $C_S$ of genus 5. One can also calculate $n_p(S)$ as
\begin{equation}
n_p(S) = \dfrac{p}{16} - \dfrac{A_{S, p \text{ mod } 4}}{16} + C_{S, p \text{ mod } 24},
\end{equation}
where $A_{S, p \text{ mod } 4}$ is the main error term written in terms of traces of Frobenius and $C_{S, p \text{ mod } 24}$ is a small correction term bounded by an absolute constant. The values of $A_{S, p \text{ mod } 4}$ and $C_{S, p \text{ mod } 24}$ are given in Tables \ref{l4table} and \ref{l4errors}, respectively.
\end{Theorem}

\begingroup
\small
\setlength{\tabcolsep}{2pt}
\renewcommand{\arraystretch}{1.5}

\begin{longtable}{|l|*{2}{c|}}
\caption{Simplified error terms $A_{S, p \text{ mod } 4}$ for $l = 4$.}
\label{l4table}\\

\hline
$S$ & $4k+1$ & $4k+3$\\
\hline
\endfirsthead

\hline
$S$ & $4k+1$ & $4k+3$\\
\hline
\endhead

\hline
\multicolumn{3}{r}{\textit{Continued on the next page}}\\
\endfoot

\hline
\endlastfoot

$\mathrm{RRRR}$ & $2 a_{0} + 2 a_{1} + a_{4}$ & $a_{4}$ \\
$\mathrm{RRNR}$ & $-2 a_{0} - a_{4}$ & $2 a_{1} - a_{4}$ \\
$\mathrm{RNRR}$ & $-2 a_{0} - a_{4}$ & $-2 a_{1} - a_{4}$ \\
$\mathrm{NRRR}$ & $-2 a_{1} - a_{4}$ & $-a_{4}$ \\
$\mathrm{RRRN}$ & $-2 a_{1} - a_{4}$ & $-a_{4}$ \\
$\mathrm{RRNN}$ & $a_{4}$ & $-2 a_{1} + a_{4}$ \\
$\mathrm{NNRR}$ & $a_{4}$ & $2 a_{1} + a_{4}$ \\
$\mathrm{NRRN}$ & $-2 a_{0} + 2 a_{1} + a_{4}$ & $a_{4}$ \\
$\mathrm{RNNR}$ & $2 a_{0} - 2 a_{1} + a_{4}$ & $a_{4}$ \\
$\mathrm{RNRN}$ & $a_{4}$ & $2 a_{1} + a_{4}$ \\
$\mathrm{NRNR}$ & $a_{4}$ & $-2 a_{1} + a_{4}$ \\
$\mathrm{NNNR}$ & $2 a_{1} - a_{4}$ & $-a_{4}$ \\
$\mathrm{RNNN}$ & $2 a_{1} - a_{4}$ & $-a_{4}$ \\
$\mathrm{NRNN}$ & $2 a_{0} - a_{4}$ & $2 a_{1} - a_{4}$ \\
$\mathrm{NNRN}$ & $2 a_{0} - a_{4}$ & $-2 a_{1} - a_{4}$ \\
$\mathrm{NNNN}$ & $-2 a_{0} - 2 a_{1} + a_{4}$ & $a_{4}$ \\

\end{longtable}
\endgroup

\begingroup
\small
\setlength{\tabcolsep}{2pt}
\renewcommand{\arraystretch}{1.6}

\begin{longtable}{|l|*{8}{c|}}
\caption{Correction terms $C_{S, p \text{ mod } 24}$ for $l = 4$.}
\label{l4errors}\\

\hline
$S$ & $24k+1$ & $24k+5$ & $24k+7$ & $24k+11$ & $24k+13$ & $24k+17$ & $24k+19$ & $24k + 23$\\
\hline
\endfirsthead

\hline
$S$ & $24k+1$ & $24k+5$ & $24k+7$ & $24k+11$ & $24k+13$ & $24k+17$ & $24k+19$ & $24k + 23$\\
\hline
\endhead

\hline
\multicolumn{9}{r}{\textit{Continued on the next page}}\\
\endfoot

\hline
\endlastfoot

$\mathrm{RRRR}$ & $-\frac{39}{16}$ & $-\frac{7}{16}$ & $-\frac{7}{16}$ & $-\frac{7}{16}$ & $-\frac{7}{16}$ & $-\frac{23}{16}$ & $-\frac{7}{16}$ & $-\frac{15}{16}$ \\
$\mathrm{NRRR}$ & $-\frac{7}{16}$ & $-\frac{7}{16}$ & $-\frac{7}{16}$ & $\frac{1}{16}$ & $-\frac{7}{16}$ & $-\frac{7}{16}$ & $\frac{1}{16}$ & $-\frac{15}{16}$ \\
$\mathrm{RNRR}$ & $-\frac{7}{16}$ & $\frac{1}{16}$ & $\frac{1}{16}$ & $-\frac{7}{16}$ & $-\frac{7}{16}$ & $-\frac{7}{16}$ & $-\frac{7}{16}$ & $\frac{1}{16}$ \\
$\mathrm{NNRR}$ & $\frac{1}{16}$ & $-\frac{7}{16}$ & $-\frac{15}{16}$ & $\frac{1}{16}$ & $\frac{1}{16}$ & $\frac{1}{16}$ & $\frac{1}{16}$ & $-\frac{15}{16}$ \\
$\mathrm{RRNR}$ & $-\frac{7}{16}$ & $\frac{1}{16}$ & $\frac{1}{16}$ & $-\frac{7}{16}$ & $-\frac{7}{16}$ & $-\frac{7}{16}$ & $-\frac{7}{16}$ & $\frac{1}{16}$ \\
$\mathrm{NRNR}$ & $\frac{1}{16}$ & $-\frac{7}{16}$ & $\frac{1}{16}$ & $-\frac{15}{16}$ & $-\frac{15}{16}$ & $\frac{1}{16}$ & $\frac{1}{16}$ & $\frac{1}{16}$ \\
$\mathrm{RNNR}$ & $\frac{1}{16}$ & $\frac{1}{16}$ & $-\frac{7}{16}$ & $-\frac{7}{16}$ & $\frac{1}{16}$ & $\frac{1}{16}$ & $-\frac{7}{16}$ & $\frac{1}{16}$ \\
$\mathrm{NNNR}$ & $\frac{1}{16}$ & $\frac{1}{16}$ & $-\frac{7}{16}$ & $\frac{1}{16}$ & $\frac{1}{16}$ & $\frac{1}{16}$ & $\frac{1}{16}$ & $-\frac{15}{16}$ \\
$\mathrm{RRRN}$ & $-\frac{7}{16}$ & $-\frac{7}{16}$ & $-\frac{7}{16}$ & $\frac{1}{16}$ & $-\frac{7}{16}$ & $-\frac{7}{16}$ & $\frac{1}{16}$ & $\frac{1}{16}$ \\
$\mathrm{NRRN}$ & $\frac{1}{16}$ & $\frac{1}{16}$ & $-\frac{7}{16}$ & $-\frac{7}{16}$ & $\frac{1}{16}$ & $-\frac{15}{16}$ & $-\frac{7}{16}$ & $\frac{1}{16}$ \\
$\mathrm{RNRN}$ & $\frac{1}{16}$ & $-\frac{7}{16}$ & $\frac{1}{16}$ & $\frac{1}{16}$ & $-\frac{15}{16}$ & $\frac{1}{16}$ & $\frac{1}{16}$ & $\frac{1}{16}$ \\
$\mathrm{NNRN}$ & $\frac{1}{16}$ & $-\frac{7}{16}$ & $\frac{1}{16}$ & $-\frac{7}{16}$ & $\frac{1}{16}$ & $\frac{1}{16}$ & $-\frac{7}{16}$ & $\frac{1}{16}$ \\
$\mathrm{RRNN}$ & $\frac{1}{16}$ & $-\frac{7}{16}$ & $\frac{1}{16}$ & $\frac{1}{16}$ & $\frac{1}{16}$ & $\frac{1}{16}$ & $-\frac{15}{16}$ & $\frac{1}{16}$ \\
$\mathrm{NRNN}$ & $\frac{1}{16}$ & $-\frac{7}{16}$ & $\frac{1}{16}$ & $-\frac{7}{16}$ & $\frac{1}{16}$ & $\frac{1}{16}$ & $-\frac{7}{16}$ & $\frac{1}{16}$ \\
$\mathrm{RNNN}$ & $\frac{1}{16}$ & $\frac{1}{16}$ & $-\frac{7}{16}$ & $\frac{1}{16}$ & $\frac{1}{16}$ & $\frac{1}{16}$ & $\frac{1}{16}$ & $\frac{1}{16}$ \\
$\mathrm{NNNN}$ & $-\frac{7}{16}$ & $-\frac{7}{16}$ & $-\frac{7}{16}$ & $-\frac{7}{16}$ & $-\frac{7}{16}$ & $-\frac{7}{16}$ & $-\frac{7}{16}$ & $-\frac{15}{16}$ \\

\end{longtable}
\endgroup

\begin{Example}
Let $p = 24k+17$, $S = (RNRN)$. Then we have
\begin{equation}
n_p(S) = \dfrac{p}{16} - \dfrac{A_{RNRN, 1}}{16} + C_{RNRN, 17} = \dfrac{p-a_4+1}{16}.
\end{equation}
\end{Example}

\section{Case $l = 5$}
For $l = 5$ we have several more curves with nonzero trace of Frobenius:
\begin{equation}
E_5: y^2 = x(x+1)(x+4); E_6: y^2 = x(x+2)(x+4);
\end{equation}
\begin{equation}
E_7: y^2 = x(x+3)(x+4); E_8: y^2 = (x+1)(x+2)(x+4);
\end{equation}
\begin{equation}
E_9: y^2 = (x+1)(x+3)(x+4); E_{10}: y^2 = (x+2)(x+3)(x+4);
\end{equation}
\begin{equation}
E_{11}: y^2 = x(x+1)(x+2)(x+4); E_{12}: y^2 = x(x+1)(x+3)(x+4);
\end{equation}
\begin{equation}
E_{13}: y^2 = x(x+2)(x+3)(x+4); E_{14}: y^2 = (x+1)(x+2)(x+3)(x+4);
\end{equation}
\begin{equation}
C: y^2 = x(x+1)(x+2)(x+3)(x+4).
\end{equation}

A direct computation gives $E_0 \cong E_{10}$, $E_2 \cong E_9$, $E_1 \cong E_8$ over $\mathbb{Q}$. Curves $E_7$ and $E_5$ are isomorphic over $\mathbb{F}_p$ if $p = 4k + 1$ and one is the quadratic twist of the other if $p = 4k+3$; $E_{6}$ is isomorphic to its own quadratic twist if $p = 4k+3$ and isomorphic to $E_0$ if $2$ is a quadratic residue modulo $p$, and to the quadratic twist of $E_0$ otherwise. It is also easy to see that $E_4 \cong E_{14}$ and $E_{11} \cong E_{13}$ over $\mathbb{Q}$. Therefore, we have:
\begin{equation}
a_{10} = a_0, a_9 = a_2, a_8 = a_1, a_{14} = a_4, a_{13} = a_{11},
\end{equation}
\begin{equation}
a_7 = (-1)^{\frac{p-1}{2}}a_5, a_{6} = (-1)^{\frac{p-1}{2}} a_{6} = (-1)^{\frac{p^2-1}{8}} a_0.
\end{equation}

We will also show that the trace $a_C$ of the curve $C$ may be nonzero only if $p = 8k+1$, which is a little trickier. If $p = 4k + 3$, then the linear substitution $x \mapsto -(x + 4), y \mapsto y$ shows that $C$ is isomorphic to its quadratic twist over $\mathbb{F}_p$ and therefore $a_{C} = 0$. Now, if we shift the curve by $x \mapsto x - 2$, looking at its smooth projective closure given in the weighted projective plane $\mathbb{P}(1, 3, 1)$ by the equation
\begin{equation}
Y^2 = Z(X-2Z)(X-Z)X(X+Z)(X+2Z)
\end{equation}
and considering it in the affine map $\{X \ne 0\}$, we see that it is given by the equation
\begin{equation}
y^2 = z(1-2z)(1-z)(1+z)(1+2z),
\end{equation}
which is equivalent to
\begin{equation}
8y^2 = 2z(1-2z)(2-2z)(2+2z)(1+2z).
\end{equation}
Now, taking $z \mapsto \frac{z}{2}, y \mapsto \frac{y}{2}$ yields that $C$ is isomorphic to its quadratic twist if $2$ is a quadratic nonresidue modulo $p$. Thus, if $2$ is a quadratic nonresidue modulo $p$, that is, if $p \equiv 3$ or $5\;(\mbox{mod}\;8)$, then $a_{C} = 0$. Following \cite[Section 4.3]{vladuc}, we note that the given composition of transformations is a non-hyperelliptic involution of $C$, which yields an isogeny over $\mathbb{Q}(\sqrt{2})$ from the Jacobian of $C$ to the product $\widetilde{E_{15}} \times \widetilde{E_{16}}$ of two elliptic curves $\widetilde{E_{15}}$ and $\widetilde{E_{16}}$ given by the equations
\begin{equation}
\widetilde{E_{15}}: y^2 = x^3 + 2\sqrt{2}x^2 - 9x - 18\sqrt{2},
\end{equation}
\begin{equation}
\widetilde{E_{16}}: y^2 = x^3 - 2\sqrt{2}x^2 - 9x + 18\sqrt{2}.
\end{equation}
Consider the curve
\begin{equation}
E_{15}: y^2 = x^3+2x^2-9x/2-9.
\end{equation}
Notice that $\widetilde{E_{15}}$ is isomorphic to the $\sqrt{2}$-twist of $E_{15}$ and $\widetilde{E_{16}}$ is isomorphic to the $(-\sqrt{2})$-twist of $E_{15}$. Therefore, we have
\begin{equation}
a_C = \tilde{a_{15}} + \tilde{a_{16}} = \left(\leg{\sqrt{2}} + \leg{-\sqrt{2}}\right)a_{15},
\end{equation}
which equals to $2\left(\frac{\sqrt{2}}{p}\right)a_{15}$ if $p = 8k+1$.
\begin{Remark}
    Our curve $E_{15}$ is denoted by $E_{16}$ in \cite[Section 4.3]{vladuc}.
\end{Remark}
 
Similar transformations give
\begin{equation}
    E_5 \cong_\mathbb{Q} E_4,\quad E_{11} \cong_{\mathbb{Q}(\sqrt{2})} E_1,
\end{equation}
which leads to $a_5 = a_4$ and $a_{11} = (-1)^\frac{p^2-1}{8}a_1$.

As in the previous section, we summarize these relations in the following theorem.
\begin{Theorem} For a prime $p \ge 5$ and a pattern $S$ of consecutive quadratic residues and quadratic nonresidues of length $l = 5$, the number $n_p(S)$ of occurrences of $S$ in $W_p$ is controlled (in the sense of equation (\ref{number_of_points_and_np})) by the curve $C_S$ of genus 17. We have
\begin{equation}
n_p(S) = \dfrac{p}{32} - \dfrac{A_{S, p \text{ mod } 8}}{32} + C_{S, p \text{ mod } 24},
\end{equation}
where $A_{S, p \text{ mod } 8}$ is the error term written in terms of traces of Frobenius and $C_{S, p \text{ mod } 24}$ is a small correction term. The values of $A_{S, p \text{ mod } 8}$ and $C_{S, p \text{ mod } 24}$ are given in Tables \ref{l5table} and \ref{l5errors}, respectively.
\end{Theorem}

\begingroup
\small
\setlength{\tabcolsep}{2pt}
\renewcommand{\arraystretch}{1.6}

\begin{longtable}{|l|*{4}{c|}}
\caption{Simplified error terms $A_{S, p \text{ mod } 8}$ for $l = 5$, $s = 2\left(\frac{\sqrt{2}}{p}\right)$.}
\label{l5table}\\

\hline
$S$ & $8k+1$ & $8k+3$ & $8k+5$ & $8k+7$ \\
\hline
\endfirsthead

\hline
$S$ & $8k+1$ & $8k+3$ & $8k+5$ & $8k+7$ \\
\hline
\endhead

\hline
\multicolumn{5}{r}{\textit{Continued on the next page}}\\
\endfoot

\hline
\endlastfoot

$\mathrm{RRRRR}$ & $4 a_{0} + 6 a_{1} + 4 a_{4} + a_{12} + sa_{15}$ & $-2 a_{1} + 2 a_{4} + a_{12}$ & $2 a_{0} + 2 a_{1} + 4 a_{4} + a_{12}$ & $2 a_{1} + 2 a_{4} + a_{12}$ \\
$\mathrm{NRRRR}$ & $-2 a_{1} - 2 a_{4} - a_{12} - sa_{15}$ & $2 a_{1} - a_{12}$ & $2 a_{0} + 2 a_{1} - 2 a_{4} - a_{12}$ & $-2 a_{1} - a_{12}$ \\
$\mathrm{RNRRR}$ & $-2 a_{1} - 2 a_{4} - a_{12} - sa_{15}$ & $-2 a_{1} - 4 a_{4} - a_{12}$ & $-2 a_{0} - 2 a_{1} - 2 a_{4} - a_{12}$ & $-2 a_{1} - 4 a_{4} - a_{12}$ \\
$\mathrm{NNRRR}$ & $-2 a_{1} + a_{12} + sa_{15}$ & $2 a_{1} + 2 a_{4} + a_{12}$ & $2 a_{0} - 2 a_{1} + a_{12}$ & $2 a_{1} + 2 a_{4} + a_{12}$ \\
$\mathrm{RRNRR}$ & $-4 a_{0} - 2 a_{1} + a_{12} - sa_{15}$ & $2 a_{1} - 2 a_{4} + a_{12}$ & $-2 a_{0} + 2 a_{1} + a_{12}$ & $-2 a_{1} - 2 a_{4} + a_{12}$ \\
$\mathrm{NRNRR}$ & $2 a_{1} - 2 a_{4} - a_{12} + sa_{15}$ & $-6 a_{1} - a_{12}$ & $-2 a_{0} - 2 a_{1} - 2 a_{4} - a_{12}$ & $-2 a_{1} - a_{12}$ \\
$\mathrm{RNNRR}$ & $-2 a_{1} + 2 a_{4} - a_{12} + sa_{15}$ & $2 a_{1} - a_{12}$ & $2 a_{0} - 2 a_{1} + 2 a_{4} - a_{12}$ & $2 a_{1} - a_{12}$ \\
$\mathrm{NNNRR}$ & $2 a_{1} + a_{12} - sa_{15}$ & $2 a_{1} + 2 a_{4} + a_{12}$ & $-2 a_{0} + 2 a_{1} + a_{12}$ & $2 a_{1} + 2 a_{4} + a_{12}$ \\
$\mathrm{RRRNR}$ & $-2 a_{1} - 2 a_{4} - a_{12} - sa_{15}$ & $2 a_{1} - a_{12}$ & $-2 a_{0} - 2 a_{1} - 2 a_{4} - a_{12}$ & $2 a_{1} - a_{12}$ \\
$\mathrm{NRRNR}$ & $-4 a_{0} + 2 a_{1} + a_{12} + sa_{15}$ & $2 a_{1} - 2 a_{4} + a_{12}$ & $-2 a_{0} + 2 a_{1} + a_{12}$ & $2 a_{1} - 2 a_{4} + a_{12}$ \\
$\mathrm{RNRNR}$ & $-2 a_{1} + a_{12} + sa_{15}$ & $2 a_{1} + 2 a_{4} + a_{12}$ & $-2 a_{0} + 2 a_{1} + a_{12}$ & $-2 a_{1} + 2 a_{4} + a_{12}$ \\
$\mathrm{NNRNR}$ & $2 a_{1} + 2 a_{4} - a_{12} - sa_{15}$ & $-6 a_{1} - a_{12}$ & $2 a_{0} - 2 a_{1} + 2 a_{4} - a_{12}$ & $-2 a_{1} - a_{12}$ \\
$\mathrm{RRNNR}$ & $-2 a_{1} + 2 a_{4} - a_{12} + sa_{15}$ & $-2 a_{1} + 4 a_{4} - a_{12}$ & $2 a_{0} - 2 a_{1} + 2 a_{4} - a_{12}$ & $-2 a_{1} + 4 a_{4} - a_{12}$ \\
$\mathrm{NRNNR}$ & $4 a_{0} - 2 a_{1} + a_{12} - sa_{15}$ & $2 a_{1} - 2 a_{4} + a_{12}$ & $2 a_{0} - 2 a_{1} + a_{12}$ & $2 a_{1} - 2 a_{4} + a_{12}$ \\
$\mathrm{RNNNR}$ & $6 a_{1} - 4 a_{4} + a_{12} - sa_{15}$ & $-2 a_{1} - 2 a_{4} + a_{12}$ & $2 a_{0} + 2 a_{1} - 4 a_{4} + a_{12}$ & $2 a_{1} - 2 a_{4} + a_{12}$ \\
$\mathrm{NNNNR}$ & $-2 a_{1} + 2 a_{4} - a_{12} + sa_{15}$ & $2 a_{1} - a_{12}$ & $-2 a_{0} + 2 a_{1} + 2 a_{4} - a_{12}$ & $-2 a_{1} - a_{12}$ \\
$\mathrm{RRRRN}$ & $-2 a_{1} - 2 a_{4} - a_{12} - sa_{15}$ & $2 a_{1} - a_{12}$ & $2 a_{0} + 2 a_{1} - 2 a_{4} - a_{12}$ & $-2 a_{1} - a_{12}$ \\
$\mathrm{NRRRN}$ & $-2 a_{1} + a_{12} + sa_{15}$ & $-2 a_{1} - 2 a_{4} + a_{12}$ & $-2 a_{0} - 6 a_{1} + a_{12}$ & $2 a_{1} - 2 a_{4} + a_{12}$ \\
$\mathrm{RNRRN}$ & $-4 a_{0} + 2 a_{1} + a_{12} + sa_{15}$ & $-2 a_{1} + 2 a_{4} + a_{12}$ & $-2 a_{0} + 2 a_{1} + a_{12}$ & $-2 a_{1} + 2 a_{4} + a_{12}$ \\
$\mathrm{NNRRN}$ & $2 a_{1} + 2 a_{4} - a_{12} - sa_{15}$ & $2 a_{1} - a_{12}$ & $-2 a_{0} + 2 a_{1} + 2 a_{4} - a_{12}$ & $2 a_{1} - a_{12}$ \\
$\mathrm{RRNRN}$ & $2 a_{1} - 2 a_{4} - a_{12} + sa_{15}$ & $2 a_{1} - a_{12}$ & $-2 a_{0} - 2 a_{1} - 2 a_{4} - a_{12}$ & $6 a_{1} - a_{12}$ \\
$\mathrm{NRNRN}$ & $-2 a_{1} + 4 a_{4} + a_{12} - sa_{15}$ & $2 a_{1} + 2 a_{4} + a_{12}$ & $2 a_{0} + 2 a_{1} + 4 a_{4} + a_{12}$ & $-2 a_{1} + 2 a_{4} + a_{12}$ \\
$\mathrm{RNNRN}$ & $4 a_{0} - 2 a_{1} + a_{12} - sa_{15}$ & $-2 a_{1} + 2 a_{4} + a_{12}$ & $2 a_{0} - 2 a_{1} + a_{12}$ & $-2 a_{1} + 2 a_{4} + a_{12}$ \\
$\mathrm{NNNRN}$ & $2 a_{1} - 2 a_{4} - a_{12} + sa_{15}$ & $-2 a_{1} - 4 a_{4} - a_{12}$ & $2 a_{0} + 2 a_{1} - 2 a_{4} - a_{12}$ & $-2 a_{1} - 4 a_{4} - a_{12}$ \\
$\mathrm{RRRNN}$ & $-2 a_{1} + a_{12} + sa_{15}$ & $-2 a_{1} - 2 a_{4} + a_{12}$ & $2 a_{0} - 2 a_{1} + a_{12}$ & $-2 a_{1} - 2 a_{4} + a_{12}$ \\
$\mathrm{NRRNN}$ & $2 a_{1} + 2 a_{4} - a_{12} - sa_{15}$ & $-2 a_{1} + 4 a_{4} - a_{12}$ & $-2 a_{0} + 2 a_{1} + 2 a_{4} - a_{12}$ & $-2 a_{1} + 4 a_{4} - a_{12}$ \\
$\mathrm{RNRNN}$ & $2 a_{1} + 2 a_{4} - a_{12} - sa_{15}$ & $2 a_{1} - a_{12}$ & $2 a_{0} - 2 a_{1} + 2 a_{4} - a_{12}$ & $6 a_{1} - a_{12}$ \\
$\mathrm{NNRNN}$ & $4 a_{0} - 2 a_{1} - 4 a_{4} + a_{12} + sa_{15}$ & $2 a_{1} - 2 a_{4} + a_{12}$ & $2 a_{0} + 2 a_{1} - 4 a_{4} + a_{12}$ & $-2 a_{1} - 2 a_{4} + a_{12}$ \\
$\mathrm{RRNNN}$ & $2 a_{1} + a_{12} - sa_{15}$ & $-2 a_{1} - 2 a_{4} + a_{12}$ & $-2 a_{0} + 2 a_{1} + a_{12}$ & $-2 a_{1} - 2 a_{4} + a_{12}$ \\
$\mathrm{NRNNN}$ & $2 a_{1} - 2 a_{4} - a_{12} + sa_{15}$ & $2 a_{1} - a_{12}$ & $2 a_{0} + 2 a_{1} - 2 a_{4} - a_{12}$ & $2 a_{1} - a_{12}$ \\
$\mathrm{RNNNN}$ & $-2 a_{1} + 2 a_{4} - a_{12} + sa_{15}$ & $2 a_{1} - a_{12}$ & $-2 a_{0} + 2 a_{1} + 2 a_{4} - a_{12}$ & $-2 a_{1} - a_{12}$ \\
$\mathrm{NNNNN}$ & $-4 a_{0} - 2 a_{1} + a_{12} - sa_{15}$ & $-2 a_{1} + 2 a_{4} + a_{12}$ & $-2 a_{0} - 6 a_{1} + a_{12}$ & $2 a_{1} + 2 a_{4} + a_{12}$ \\

\end{longtable}
\endgroup

\begingroup
\small
\setlength{\tabcolsep}{2pt}
\renewcommand{\arraystretch}{1.6}

\begin{longtable}{|l|*{8}{c|}}
\caption{Correction terms $C_{S, p \text{ mod } 24}$ for $l = 5$.}
\label{l5errors}\\

\hline
$S$ & $24k+1$ & $24k+5$ & $24k+7$ & $24k+11$ & $24k+13$ & $24k+17$ & $24k+19$ & $24k + 23$\\
\hline
\endfirsthead

\hline
$S$ & $24k+1$ & $24k+5$ & $24k+7$ & $24k+11$ & $24k+13$ & $24k+17$ & $24k+19$ & $24k + 23$\\
\hline
\endhead

\hline
\multicolumn{9}{r}{\textit{Continued on the next page}}\\
\endfoot

\hline
\endlastfoot

$\mathrm{RRRRR}$ & $-\frac{95}{32}$ & $-\frac{15}{32}$ & $-\frac{15}{32}$ & $-\frac{15}{32}$ & $-\frac{15}{32}$ & $-\frac{31}{32}$ & $-\frac{15}{32}$ & $-\frac{31}{32}$ \\
$\mathrm{NRRRR}$ & $-\frac{15}{32}$ & $\frac{1}{32}$ & $\frac{1}{32}$ & $\frac{1}{32}$ & $\frac{1}{32}$ & $-\frac{15}{32}$ & $\frac{1}{32}$ & $-\frac{31}{32}$ \\
$\mathrm{RNRRR}$ & $-\frac{15}{32}$ & $\frac{1}{32}$ & $\frac{1}{32}$ & $\frac{1}{32}$ & $-\frac{15}{32}$ & $-\frac{15}{32}$ & $\frac{1}{32}$ & $\frac{1}{32}$ \\
$\mathrm{NNRRR}$ & $\frac{1}{32}$ & $-\frac{15}{32}$ & $-\frac{15}{32}$ & $\frac{1}{32}$ & $\frac{1}{32}$ & $\frac{1}{32}$ & $\frac{1}{32}$ & $-\frac{31}{32}$ \\
$\mathrm{RRNRR}$ & $-\frac{15}{32}$ & $\frac{1}{32}$ & $\frac{1}{32}$ & $-\frac{15}{32}$ & $-\frac{31}{32}$ & $-\frac{15}{32}$ & $-\frac{15}{32}$ & $\frac{1}{32}$ \\
$\mathrm{NRNRR}$ & $\frac{1}{32}$ & $\frac{1}{32}$ & $\frac{1}{32}$ & $-\frac{31}{32}$ & $-\frac{15}{32}$ & $\frac{1}{32}$ & $\frac{1}{32}$ & $\frac{1}{32}$ \\
$\mathrm{RNNRR}$ & $\frac{1}{32}$ & $-\frac{15}{32}$ & $-\frac{15}{32}$ & $\frac{1}{32}$ & $\frac{1}{32}$ & $\frac{1}{32}$ & $\frac{1}{32}$ & $\frac{1}{32}$ \\
$\mathrm{NNNRR}$ & $\frac{1}{32}$ & $\frac{1}{32}$ & $-\frac{15}{32}$ & $\frac{1}{32}$ & $\frac{1}{32}$ & $\frac{1}{32}$ & $\frac{1}{32}$ & $-\frac{31}{32}$ \\
$\mathrm{RRRNR}$ & $-\frac{15}{32}$ & $\frac{1}{32}$ & $-\frac{15}{32}$ & $\frac{1}{32}$ & $-\frac{15}{32}$ & $-\frac{15}{32}$ & $\frac{1}{32}$ & $\frac{1}{32}$ \\
$\mathrm{NRRNR}$ & $\frac{1}{32}$ & $\frac{1}{32}$ & $-\frac{15}{32}$ & $-\frac{15}{32}$ & $\frac{1}{32}$ & $-\frac{31}{32}$ & $-\frac{15}{32}$ & $\frac{1}{32}$ \\
$\mathrm{RNRNR}$ & $\frac{1}{32}$ & $\frac{1}{32}$ & $\frac{1}{32}$ & $\frac{1}{32}$ & $-\frac{31}{32}$ & $\frac{1}{32}$ & $\frac{1}{32}$ & $\frac{1}{32}$ \\
$\mathrm{NNRNR}$ & $\frac{1}{32}$ & $-\frac{15}{32}$ & $\frac{1}{32}$ & $-\frac{31}{32}$ & $\frac{1}{32}$ & $\frac{1}{32}$ & $\frac{1}{32}$ & $\frac{1}{32}$ \\
$\mathrm{RRNNR}$ & $\frac{1}{32}$ & $-\frac{15}{32}$ & $\frac{1}{32}$ & $\frac{1}{32}$ & $\frac{1}{32}$ & $\frac{1}{32}$ & $-\frac{31}{32}$ & $\frac{1}{32}$ \\
$\mathrm{NRNNR}$ & $\frac{1}{32}$ & $-\frac{15}{32}$ & $-\frac{15}{32}$ & $-\frac{15}{32}$ & $\frac{1}{32}$ & $\frac{1}{32}$ & $-\frac{15}{32}$ & $\frac{1}{32}$ \\
$\mathrm{RNNNR}$ & $\frac{1}{32}$ & $\frac{1}{32}$ & $-\frac{15}{32}$ & $\frac{1}{32}$ & $\frac{1}{32}$ & $\frac{1}{32}$ & $\frac{1}{32}$ & $\frac{1}{32}$ \\
$\mathrm{NNNNR}$ & $\frac{1}{32}$ & $\frac{1}{32}$ & $\frac{1}{32}$ & $\frac{1}{32}$ & $\frac{1}{32}$ & $\frac{1}{32}$ & $\frac{1}{32}$ & $-\frac{31}{32}$ \\
$\mathrm{RRRRN}$ & $-\frac{15}{32}$ & $\frac{1}{32}$ & $\frac{1}{32}$ & $\frac{1}{32}$ & $\frac{1}{32}$ & $-\frac{15}{32}$ & $\frac{1}{32}$ & $\frac{1}{32}$ \\
$\mathrm{NRRRN}$ & $\frac{1}{32}$ & $-\frac{15}{32}$ & $-\frac{15}{32}$ & $\frac{1}{32}$ & $-\frac{15}{32}$ & $\frac{1}{32}$ & $\frac{1}{32}$ & $\frac{1}{32}$ \\
$\mathrm{RNRRN}$ & $\frac{1}{32}$ & $\frac{1}{32}$ & $\frac{1}{32}$ & $-\frac{15}{32}$ & $\frac{1}{32}$ & $-\frac{31}{32}$ & $-\frac{15}{32}$ & $\frac{1}{32}$ \\
$\mathrm{NNRRN}$ & $\frac{1}{32}$ & $\frac{1}{32}$ & $-\frac{15}{32}$ & $\frac{1}{32}$ & $\frac{1}{32}$ & $\frac{1}{32}$ & $\frac{1}{32}$ & $\frac{1}{32}$ \\
$\mathrm{RRNRN}$ & $\frac{1}{32}$ & $\frac{1}{32}$ & $\frac{1}{32}$ & $\frac{1}{32}$ & $-\frac{15}{32}$ & $\frac{1}{32}$ & $\frac{1}{32}$ & $\frac{1}{32}$ \\
$\mathrm{NRNRN}$ & $\frac{1}{32}$ & $-\frac{15}{32}$ & $\frac{1}{32}$ & $\frac{1}{32}$ & $-\frac{15}{32}$ & $\frac{1}{32}$ & $\frac{1}{32}$ & $\frac{1}{32}$ \\
$\mathrm{RNNRN}$ & $\frac{1}{32}$ & $-\frac{15}{32}$ & $\frac{1}{32}$ & $-\frac{15}{32}$ & $\frac{1}{32}$ & $\frac{1}{32}$ & $-\frac{15}{32}$ & $\frac{1}{32}$ \\
$\mathrm{NNNRN}$ & $\frac{1}{32}$ & $\frac{1}{32}$ & $\frac{1}{32}$ & $\frac{1}{32}$ & $\frac{1}{32}$ & $\frac{1}{32}$ & $\frac{1}{32}$ & $\frac{1}{32}$ \\
$\mathrm{RRRNN}$ & $\frac{1}{32}$ & $-\frac{15}{32}$ & $\frac{1}{32}$ & $\frac{1}{32}$ & $\frac{1}{32}$ & $\frac{1}{32}$ & $\frac{1}{32}$ & $\frac{1}{32}$ \\
$\mathrm{NRRNN}$ & $\frac{1}{32}$ & $\frac{1}{32}$ & $\frac{1}{32}$ & $\frac{1}{32}$ & $\frac{1}{32}$ & $\frac{1}{32}$ & $-\frac{31}{32}$ & $\frac{1}{32}$ \\
$\mathrm{RNRNN}$ & $\frac{1}{32}$ & $-\frac{15}{32}$ & $\frac{1}{32}$ & $\frac{1}{32}$ & $\frac{1}{32}$ & $\frac{1}{32}$ & $\frac{1}{32}$ & $\frac{1}{32}$ \\
$\mathrm{NNRNN}$ & $\frac{1}{32}$ & $\frac{1}{32}$ & $\frac{1}{32}$ & $-\frac{15}{32}$ & $\frac{1}{32}$ & $\frac{1}{32}$ & $-\frac{15}{32}$ & $\frac{1}{32}$ \\
$\mathrm{RRNNN}$ & $\frac{1}{32}$ & $\frac{1}{32}$ & $\frac{1}{32}$ & $\frac{1}{32}$ & $\frac{1}{32}$ & $\frac{1}{32}$ & $\frac{1}{32}$ & $\frac{1}{32}$ \\
$\mathrm{NRNNN}$ & $\frac{1}{32}$ & $\frac{1}{32}$ & $-\frac{15}{32}$ & $\frac{1}{32}$ & $\frac{1}{32}$ & $\frac{1}{32}$ & $\frac{1}{32}$ & $\frac{1}{32}$ \\
$\mathrm{RNNNN}$ & $\frac{1}{32}$ & $\frac{1}{32}$ & $\frac{1}{32}$ & $\frac{1}{32}$ & $\frac{1}{32}$ & $\frac{1}{32}$ & $\frac{1}{32}$ & $\frac{1}{32}$ \\
$\mathrm{NNNNN}$ & $-\frac{15}{32}$ & $-\frac{15}{32}$ & $-\frac{15}{32}$ & $-\frac{15}{32}$ & $-\frac{15}{32}$ & $-\frac{15}{32}$ & $-\frac{15}{32}$ & $-\frac{31}{32}$ \\

\end{longtable}
\endgroup

\begin{Example}
Let $p = 24k+11$, $S = (RNNRR)$. Then we have
\begin{equation}
n_p(RNNRR) = \dfrac{p}{32} - \dfrac{A_{RNNRR, 11}}{32} + C_{RNNRR, 11} = \dfrac{p - 2a_1 + a_{12} + 1}{32}.
\end{equation}
\end{Example}

\section{Statistical properties}

This section describes statistical properties of $n_p(S)$. First, recall that Weil bound, which follows from the Riemann hypothesis for curves over finite fields, gives
\begin{equation}
    |a_i(p)| \le 2g_i \sqrt{p},
\end{equation}
where $g_i$ is the genus of the curve $E_i$ with the Frobenius trace $a_i(p)$. Therefore, the normalized trace $\frac{a_i(p)}{2\sqrt{p}}$ lies in $[-g_i, g_i]$ and one can ask how it is distributed on this interval as $p$ varies.

Let $Y$ be a compact metric space equipped with a Borel probability measure $\mu$. A sequence $\{y_p\} \subset Y$, indexed by primes, is said to be \textit{equidistributed} with respect to $\mu$ if
\begin{equation}
\lim_{n \to \infty}\dfrac{\#\{p \le n\ :\  y_p \in X\}}{\#\{p \le n\}} = \mu(X)
\end{equation}
for every $\mu$-measurable set $X \subset Y$ with $\mu(\partial X) = 0$.

Fix some positive integer $d$. Suppose for every $r$ with $1 \le r \le d-1$, $(r, d) = 1$, there exists a measure $\mu_r$ on $Y$ such that
\begin{equation}
\lim_{n \to \infty}\dfrac{\#\{p \le n\ :\  p = md+r, m \in \mathbb{Z}_+, y_p \in X\}}{\#\{p \le n\ :\  p = md+r, m \in \mathbb{Z}_+\}} = \mu_r(X).
\end{equation}
Then Dirichlet's theorem on arithmetic progressions implies that $\{y_p\}$ (indexed by all primes) is equidistributed with respect to the measure $\mu$ given by
\begin{equation}
\mu = \dfrac{1}{\varphi(d)} \sum_{1 \le r \le d-1, (r,d) = 1} \mu_r.
\end{equation}

For elliptic curves the distribution of the normalized trace is known; it depends on the endomorphisms group of the curve. If a curve $E / \mathbb{Q}$ is a CM-curve, i.e. $\End_\mathbb{C}(E) \ne \mathbb{Z}$, then the sequence of normalized traces $\left\{\frac{a_E(p)}{2\sqrt{p}}\right\}$, indexed by primes of good reduction for $E$, is equidistributed on $[-1, 1]$ with respect to the measure
\begin{equation}
\lambda_{cm} = \frac{\delta_0}{2} + \frac{\mu_{cm}}{2},
\end{equation}
where $\delta_0$ is the Dirac measure at $0$ and $\mu_{cm}$ is the measure
\begin{equation}
\mu_{cm} = \dfrac{dt}{\pi\sqrt{1-t^2}}.
\end{equation}
Moreover, the term $\delta_0$ corresponds to the subsequence $\{p = 4k+3\}$ and the term $\mu_{cm}$ corresponds to the subsequence $\{p = 4k+1\}$, i.e. $\mu_1 = \mu_{cm}$ and $\mu_3 = \delta_0$ in terms of the previous paragraph.

If the curve $E$ has no complex multiplication, i.e. $\End_\mathbb{C}(E) = \mathbb{Z}$, then the sequence $\left\{\frac{a_E(p)}{2\sqrt{p}}\right\}$ is equidistributed on $[-1, 1]$ with respect to the measure
\begin{equation}
\mu_{ST} = \dfrac{2\sqrt{1-t^2}}{\pi}dt.
\end{equation}

Note that the result for curves with no complex multiplication is the original Sato-Tate conjecture. For details on these distributions, we refer the reader to \cite{sutherland_distributions}.

However, in the formulas for $n_p(S)$ the sums of the Frobenius traces of different curves occur and it is therefore natural to consider the joint distribution of the traces.

Let us recall the definition of \textit{stochastically independent} sequences. Let $Y$ be a compact measurable space and $\mu_1, \ldots, \mu_n$ be probability measures on $Y$. Let $s_i = \{a_p^i\} \subset Y$ be a sequence indexed by primes for $i = 1, \ldots, n$ and assume that $s_i$ is equidistributed with respect to $\mu_i$ for all $i = 1, \ldots, n$. Then $s_1, \ldots, s_n$ are called stochastically independent if the sequence
\begin{equation}
\left\{(a_p^1, \ldots, a_p^n)\right\} \subset Y^n
\end{equation}
is equidistributed with respect to the product measure $\mu_1 \times \ldots \times \mu_n$ on $Y^n$. 

The hypothesis is that the sequences of normalized traces of pairwise non-isogenous curves are stochastically independent. This can be deduced from the generalized Sato-Tate conjecture, see \cite{vladuc} for details. This enables us to express the distribution of a sum of normalized traces as the convolution of the measures of the individual terms. Note that since the correction term 
\begin{equation}
2^{-l}\sum_{\varnothing \ne T \subset [1, l]} (1 - r_{S, T}) - c_{p, S}(l)
\end{equation}
 is bounded by a constant depending on $l$ but not on $p$, we therefore obtain the distribution of $\frac{1}{\sqrt{p}}n_p(S) - \frac{\sqrt{p}}{2^l}$. As can be seen from Tables \ref{l4table} and \ref{l5table}, the error term depends on the traces of five curves $E_0, E_1, E_4, E_{12}, E_{15}$. These curves can be easily checked to be pairwise non-isogenous over $\overline{\mathbb{Q}}$, see \cite[Section 4.3]{vladuc} for a thorough analysis, so we will assume that GST holds for the variety
\begin{equation}
B(4) = E_0 \times E_1 \times E_4.
\end{equation}
Then, using Table \ref{l4table}, one can deduce the density function of the limiting distribution of the normalized error term. Recall that $E_0$ is a CM-curve, whereas $E_1$ and $E_4$ have no complex multiplication. This accounts for the appropriate choice of the density functions in the convolution. Note that the intervals on which the resulting distributions are supported can be easily deduced since both $\lambda_{cm}$ and $\mu_{ST}$ are supported on $[-1, 1]$ and the support of the convolution of the measures $\lambda_1, \lambda_2$ is equal to $[a_1 + a_2, b_1+b_2]$ if $\lambda_i$ is supported on $[a_i, b_i]$.

\begin{Theorem}
Assume GST for the variety $B(4)$. Then for $r \in \{1, 3\}$ and a pattern $S$ of consecutive quadratic residues and quadratic nonresidues of length $l = 4$ the sequence
\begin{equation}
\dfrac{1}{\sqrt{p}}n_p(S) - \dfrac{\sqrt{p}}{2^l}
\end{equation}
of normalized error terms indexed by primes $p$ of the form $4k+r$, excluding the primes at which $B(4)$ has bad reduction, is equidistributed on the interval $I_{S, r}$ with respect to the measure $\mu_{S, r}$. The values of $\mu_{S, r}$ and $I_{S, r}$ are given in Tables \ref{l4dist} and \ref{l4support}, respectively.
\end{Theorem}

\begingroup
\small
\setlength{\tabcolsep}{2pt}
\renewcommand{\arraystretch}{2.0}

\begin{longtable}{|l|*{2}{c|}}
\caption{Limiting distributions $\mu_{S, r}$ for patterns of length $4$.}
\label{l4dist}\\

\hline
$S$ & $4k+1$ & $4k+3$ \\
\hline
\endfirsthead

\hline
$S$ & $4k+1$ & $4k+3$ \\
\hline
\endhead

\hline
\multicolumn{3}{r}{\textit{Continued on the next page}}\\
\endfoot

\hline
\endlastfoot

$\mathrm{RRRR}$ & $\displaystyle \frac{ 512 }{ \pi^{3} } \frac{1}{\sqrt{ 1 - 16 t_{0} ^2 }} \ast \sqrt{ 1 - 16 t_{1} ^2 } \ast \sqrt{ 1 - 64 t_{4} ^2 } \,dt $ & $\displaystyle \frac{ 16 }{ \pi } \sqrt{ 1 - 64 t_{4} ^2 } \,dt $ \\
$\mathrm{NRRR}$ & $\displaystyle \frac{ 128 }{ \pi^{2} } \sqrt{ 1 - 16 t_{1} ^2 } \ast \sqrt{ 1 - 64 t_{4} ^2 } \,dt $ & $\displaystyle \frac{ 16 }{ \pi } \sqrt{ 1 - 64 t_{4} ^2 } \,dt $ \\
$\mathrm{RNRR}$ & $\displaystyle \frac{ 64 }{ \pi^{2} } \frac{1}{\sqrt{ 1 - 16 t_{0} ^2 }} \ast \sqrt{ 1 - 64 t_{4} ^2 } \,dt $ & $\displaystyle \frac{ 128 }{ \pi^{2} } \sqrt{ 1 - 16 t_{1} ^2 } \ast \sqrt{ 1 - 64 t_{4} ^2 } \,dt $ \\
$\mathrm{NNRR}$ & $\displaystyle \frac{ 16 }{ \pi } \sqrt{ 1 - 64 t_{4} ^2 } \,dt $ & $\displaystyle \frac{ 128 }{ \pi^{2} } \sqrt{ 1 - 16 t_{1} ^2 } \ast \sqrt{ 1 - 64 t_{4} ^2 } \,dt $ \\
$\mathrm{RRNR}$ & $\displaystyle \frac{ 64 }{ \pi^{2} } \frac{1}{\sqrt{ 1 - 16 t_{0} ^2 }} \ast \sqrt{ 1 - 64 t_{4} ^2 } \,dt $ & $\displaystyle \frac{ 128 }{ \pi^{2} } \sqrt{ 1 - 16 t_{1} ^2 } \ast \sqrt{ 1 - 64 t_{4} ^2 } \,dt $ \\
$\mathrm{NRNR}$ & $\displaystyle \frac{ 16 }{ \pi } \sqrt{ 1 - 64 t_{4} ^2 } \,dt $ & $\displaystyle \frac{ 128 }{ \pi^{2} } \sqrt{ 1 - 16 t_{1} ^2 } \ast \sqrt{ 1 - 64 t_{4} ^2 } \,dt $ \\
$\mathrm{RNNR}$ & $\displaystyle \frac{ 512 }{ \pi^{3} } \frac{1}{\sqrt{ 1 - 16 t_{0} ^2 }} \ast \sqrt{ 1 - 16 t_{1} ^2 } \ast \sqrt{ 1 - 64 t_{4} ^2 } \,dt $ & $\displaystyle \frac{ 16 }{ \pi } \sqrt{ 1 - 64 t_{4} ^2 } \,dt $ \\
$\mathrm{NNNR}$ & $\displaystyle \frac{ 128 }{ \pi^{2} } \sqrt{ 1 - 16 t_{1} ^2 } \ast \sqrt{ 1 - 64 t_{4} ^2 } \,dt $ & $\displaystyle \frac{ 16 }{ \pi } \sqrt{ 1 - 64 t_{4} ^2 } \,dt $ \\
$\mathrm{RRRN}$ & $\displaystyle \frac{ 128 }{ \pi^{2} } \sqrt{ 1 - 16 t_{1} ^2 } \ast \sqrt{ 1 - 64 t_{4} ^2 } \,dt $ & $\displaystyle \frac{ 16 }{ \pi } \sqrt{ 1 - 64 t_{4} ^2 } \,dt $ \\
$\mathrm{NRRN}$ & $\displaystyle \frac{ 512 }{ \pi^{3} } \frac{1}{\sqrt{ 1 - 16 t_{0} ^2 }} \ast \sqrt{ 1 - 16 t_{1} ^2 } \ast \sqrt{ 1 - 64 t_{4} ^2 } \,dt $ & $\displaystyle \frac{ 16 }{ \pi } \sqrt{ 1 - 64 t_{4} ^2 } \,dt $ \\
$\mathrm{RNRN}$ & $\displaystyle \frac{ 16 }{ \pi } \sqrt{ 1 - 64 t_{4} ^2 } \,dt $ & $\displaystyle \frac{ 128 }{ \pi^{2} } \sqrt{ 1 - 16 t_{1} ^2 } \ast \sqrt{ 1 - 64 t_{4} ^2 } \,dt $ \\
$\mathrm{NNRN}$ & $\displaystyle \frac{ 64 }{ \pi^{2} } \frac{1}{\sqrt{ 1 - 16 t_{0} ^2 }} \ast \sqrt{ 1 - 64 t_{4} ^2 } \,dt $ & $\displaystyle \frac{ 128 }{ \pi^{2} } \sqrt{ 1 - 16 t_{1} ^2 } \ast \sqrt{ 1 - 64 t_{4} ^2 } \,dt $ \\
$\mathrm{RRNN}$ & $\displaystyle \frac{ 16 }{ \pi } \sqrt{ 1 - 64 t_{4} ^2 } \,dt $ & $\displaystyle \frac{ 128 }{ \pi^{2} } \sqrt{ 1 - 16 t_{1} ^2 } \ast \sqrt{ 1 - 64 t_{4} ^2 } \,dt $ \\
$\mathrm{NRNN}$ & $\displaystyle \frac{ 64 }{ \pi^{2} } \frac{1}{\sqrt{ 1 - 16 t_{0} ^2 }} \ast \sqrt{ 1 - 64 t_{4} ^2 } \,dt $ & $\displaystyle \frac{ 128 }{ \pi^{2} } \sqrt{ 1 - 16 t_{1} ^2 } \ast \sqrt{ 1 - 64 t_{4} ^2 } \,dt $ \\
$\mathrm{RNNN}$ & $\displaystyle \frac{ 128 }{ \pi^{2} } \sqrt{ 1 - 16 t_{1} ^2 } \ast \sqrt{ 1 - 64 t_{4} ^2 } \,dt $ & $\displaystyle \frac{ 16 }{ \pi } \sqrt{ 1 - 64 t_{4} ^2 } \,dt $ \\
$\mathrm{NNNN}$ & $\displaystyle \frac{ 512 }{ \pi^{3} } \frac{1}{\sqrt{ 1 - 16 t_{0} ^2 }} \ast \sqrt{ 1 - 16 t_{1} ^2 } \ast \sqrt{ 1 - 64 t_{4} ^2 } \,dt $ & $\displaystyle \frac{ 16 }{ \pi } \sqrt{ 1 - 64 t_{4} ^2 } \,dt $ \\

\end{longtable}
\endgroup

\begingroup
\small
\setlength{\tabcolsep}{4pt}
\renewcommand{\arraystretch}{2.0}

\begin{longtable}{|l|*{2}{c|}}
\caption{Supports $I_{S, r}$ of the limiting distributions for patterns of length $4$.}
\label{l4support}\\

\hline
$S$ & $4k+1$ & $4k+3$\\
\hline
\endfirsthead

\hline
$S$ & $4k+1$ & $4k+3$\\
\hline
\endhead

\hline
\multicolumn{3}{r}{\textit{Continued on the next page}}\\
\endfoot

\hline
\endlastfoot

$\mathrm{RRRR}$ & $ \left[- \frac{5}{8} , \frac{5}{8} \right] $ & $ \left[- \frac{1}{8} , \frac{1}{8} \right] $ \\
$\mathrm{NRRR}$ & $ \left[- \frac{3}{8} , \frac{3}{8} \right] $ & $ \left[- \frac{1}{8} , \frac{1}{8} \right] $ \\
$\mathrm{RNRR}$ & $ \left[- \frac{3}{8} , \frac{3}{8} \right] $ & $ \left[- \frac{3}{8} , \frac{3}{8} \right] $ \\
$\mathrm{NNRR}$ & $ \left[- \frac{1}{8} , \frac{1}{8} \right] $ & $ \left[- \frac{3}{8} , \frac{3}{8} \right] $ \\
$\mathrm{RRNR}$ & $ \left[- \frac{3}{8} , \frac{3}{8} \right] $ & $ \left[- \frac{3}{8} , \frac{3}{8} \right] $ \\
$\mathrm{NRNR}$ & $ \left[- \frac{1}{8} , \frac{1}{8} \right] $ & $ \left[- \frac{3}{8} , \frac{3}{8} \right] $ \\
$\mathrm{RNNR}$ & $ \left[- \frac{5}{8} , \frac{5}{8} \right] $ & $ \left[- \frac{1}{8} , \frac{1}{8} \right] $ \\
$\mathrm{NNNR}$ & $ \left[- \frac{3}{8} , \frac{3}{8} \right] $ & $ \left[- \frac{1}{8} , \frac{1}{8} \right] $ \\
$\mathrm{RRRN}$ & $ \left[- \frac{3}{8} , \frac{3}{8} \right] $ & $ \left[- \frac{1}{8} , \frac{1}{8} \right] $ \\
$\mathrm{NRRN}$ & $ \left[- \frac{5}{8} , \frac{5}{8} \right] $ & $ \left[- \frac{1}{8} , \frac{1}{8} \right] $ \\
$\mathrm{RNRN}$ & $ \left[- \frac{1}{8} , \frac{1}{8} \right] $ & $ \left[- \frac{3}{8} , \frac{3}{8} \right] $ \\
$\mathrm{NNRN}$ & $ \left[- \frac{3}{8} , \frac{3}{8} \right] $ & $ \left[- \frac{3}{8} , \frac{3}{8} \right] $ \\
$\mathrm{RRNN}$ & $ \left[- \frac{1}{8} , \frac{1}{8} \right] $ & $ \left[- \frac{3}{8} , \frac{3}{8} \right] $ \\
$\mathrm{NRNN}$ & $ \left[- \frac{3}{8} , \frac{3}{8} \right] $ & $ \left[- \frac{3}{8} , \frac{3}{8} \right] $ \\
$\mathrm{RNNN}$ & $ \left[- \frac{3}{8} , \frac{3}{8} \right] $ & $ \left[- \frac{1}{8} , \frac{1}{8} \right] $ \\
$\mathrm{NNNN}$ & $ \left[- \frac{5}{8} , \frac{5}{8} \right] $ & $ \left[- \frac{1}{8} , \frac{1}{8} \right] $ \\

\end{longtable}
\endgroup

\begin{Example}
Let $S = $(RRRN). Then the sequence
\begin{equation}
\dfrac{1}{\sqrt{p}}n_p(S) - \dfrac{\sqrt{p}}{2^l},
\end{equation}
indexed by primes $p$ of the form $4k+1$, is equidistributed on the interval $I_{RRRN, 1} = \left[-\frac38, \frac38\right]$ with respect to the measure
\begin{equation}
\mu_{RRRN, 1} = \dfrac{128}{\pi^2} \sqrt{1-16t_1^2} * \sqrt{1 - 64t_4^2} dt.
\end{equation}
\end{Example}

\begin{Remark}
Assuming GST for the variety
\begin{equation}
    B(5) = E_0 \times E_1 \times E_4 \times E_{12} \times E_{15},
\end{equation}
limiting distributions for $l = 5$ can be worked out similarly, but the resulting formulas are too lengthy to present in tabular form.
\end{Remark}

\bibliographystyle{amsplain}
\bibliography{quadratic_patterns_ref}

\end{document}